%This is an AMStex file

\input amstex
\documentstyle{amsppt}

\NoRunningHeads
\magnification=1100
\baselineskip=29pt
\parskip 9pt
\pagewidth{6.0in}
\pageheight{8.0in}

\TagsOnLeft
\CenteredTagsOnSplits

\let\pro=\proclaim
\let\endpro=\endproclaim
\def\proof{{\it Proof.\ }}

\let\raw=\rightarrow
\let\parl=\partial
\let\wed=\wedge
\let\sst=\subset

\let\alp=\alpha
\let\be=\beta
\let\ga=\gamma

\let\om=\omega
\let\Om=\Omega

\let\hta=\theta

\let\tld=\tilde

\let\rarw=\rightarrow

\topmatter
\title
Special Lagrangian Tori on a Borcea-Voisin Threefold
\endtitle

\author
Peng Lu 
\endauthor

\affil
School of Mathematics\\
University of Minnesota \\
Minneapolis, MN 55455
\endaffil

%\date{}
%\enddate
\endtopmatter

\document

%%%%%%%%%%%%%%%%%%%%%%%%%%%%%%%%%%%%%%%%%%%%%%%%%%%%%%%%%%%%%%%%%%%%%%%%%%%%%%%%%%%%%%%%%%%%%%%%%%%%%%%%%%%%%%%%%%%%%%%%%%%%

%\head 0. Introduction
%\endhead

Borcea-Voisin threefolds are Calabi-Yau manifolds. They are constructed by
Borcea ({\cite B}) and Voisin ({\cite V}) in the construction of mirror manifolds. In
\cite{SYZ}, Strominger, Yau and Zaslow propose a geometric construction of mirror
manifolds using special  Lagrangian tori (called SYZ-construction below).
Using degenerate Calabi-Yau metrics Gross and Wilson show that
SYZ-construction works for any Borcea-Voisin threefolds (\cite{GW}). In this short
note we show that there are special Lagrangian tori on one family of Borcea-Voisin 
threefolds with respect to non-degenerate  Calabi-Yau metrics. These
tori cover a large part of the threefolds and are perturbations of the special
Lagrangian tori used by  Gross and Wilson in \cite{GW}. The method is studying
the degeneration of Calabi-Yau metrics using gluing (type I degeneration). There are
other examples of special Lagrangian torus in compact Calabi-Yau threefolds. 
Bryant has a beautiful construction of the special Lagrangian torus in
some quintic threefolds (see \cite{Br}). In section 2 we give a family of special
Lagrangian submanifolds which cover $K_{{\Bbb C}P^n}$.

In some sense this note is a continuation of {\cite L}. In that paper we show 
the existence of immersed special Lagrangian tori in a Kummer type threefold. In 
this paper we show the existence of embedded special Lagrangian tori in one  
family of Borcea-Voisin 
threefolds. Readers may refer to that paper for basic definitions if necessary.   

The author thanks Naichung C. Leung for helpful discussions.

\head 1.  Special Lagrangian Tori  on Borcea-Voisin Threefolds
\endhead

 Let $E_i$ be an elliptic curve with periods $1$ and $\tau_i$ ($i=1,2,3$). The 
Borcea-Voisin Threefold in discussion is the resolution of quotient $E_1
 \times E_2 \times E_3$ by ${\Bbb Z}_2 \times {\Bbb Z}_2$. We denote the
quotient by $M_0$. Let $z_1,\, z_2, \, z_3$ be 
the coordinates of the three elliptic curves. The generators of the 
two ${\Bbb Z}_2$ actions are

$$\align
&\alp: z_1 \raw -z_1 +\frac{1}{2},\,\,  z_2 \raw -z_2 +\frac{1}{2},\,\, z_3 \raw z_3, \\
& \beta:  z_1 \raw -z_1, \,\,  z_2 \raw z_2, \,\, z_3 \raw -z_3.
\endalign
$$

\noindent The fixed locus of $\alp$ is the union of sixteen elliptic curves:

$$ (\frac{2+ \tau_1 \pm 1 \pm \tau_1 }{4}, \frac{2+ \tau_2 \pm 1 \pm \tau_2}{4})
 \times E_3. $$

\noindent The fixed locus of $\beta$ is the union of sixteen elliptic curves:

$$\frac{1+ \tau_1 \pm 1 \pm \tau_1 }{4} \times E_2 \times
\frac{1+ \tau_3 \pm 1 \pm \tau_3}{4} $$

\noindent The fixed locus of $\alp \circ \beta$ is empty. Note that the 32 elliptic
curves in the fixed locus do not intersect each other and their images in $M_0$
are 16 elliptic curves. Denote the image set by $F$. 

Let  $\pi: M \rarw M_0$  be the resolution map by a single blowup of $F$. Then $M$ 
is a Borcea-Voisin threefold. To see this fact, first we resolve the quotient of 
 $E_1 \times E_2 \times E_3$ by $\alp$. According to Kummer construction, we get
 $K \times E_3$, where $K$ is a Kummer surface. Then $\beta$ induces an action on
 $K \times E_3$. This is exactly the action used in Borcea-Voisin threefold 
construction (see {\cite B} or {\cite V}). $M$ is the Borcea-Voisin threefold
 constructed from  $K \times E_3$ and $\beta$. Note that the fixed locus of the 
involution on $K$ is two tori and $M$ is self mirror.

{\bf 1a}. Holomorphic (3,0)-form on $M$. Note that $dz_1 \wed dz_2 
\wed dz_3$ is a holomorphic (3,0)-form on $E_1 \times E_2 
\times E_3$. Let $\Omega_0$ be the induced holomorphic (3,0)-form on $M_0$.
 Denote $\pi^* \Omega_0$ by $\Omega$. Then

\pro{Lemma 1}
$\Omega$ is a holomorphic (3,0)-form on $M$.
\endpro 

\proof
We only need to check that  $\pi^* \Omega_0$ extends across the exceptional divisors 
and is nonzero everywhere on the exceptional divisors. Since $M$ is resolved by a 
single blow-up along the singular elliptic curves in $M_0$. In the normal
 direction the singularities are 
of the form ${\Bbb C}^2/{\Bbb Z}_2$. We know the required extension is possible for
 ${\Bbb C}^2/{\Bbb Z}_2$ (see, for example, {\cite L}, the proof of Lemma 3.1).   
\qed

{\bf 1b}. Ricci-flat metrics on $M$. First we describe the Ricci-flat metrics on the
total spaces of the normal bundles of the 
exceptional divisors in $M$. There are two types of exceptional divisors in $M$:
8 copies of ${\Bbb C}P^1 \times E_3$ and 8 copies of ${\Bbb C}P^1 \times E_2$. Let $E$
 be either $E_2$ or $E_3$. The total space of the normal bundle of ${\Bbb C}P^1 
\times E$ in $M$ is $K_{{\Bbb C}P^1} \times E$.  We can identify $K_{{\Bbb C}P^1}
 \times E$ with  the resolution of $({\Bbb C}^2/{\Bbb Z}_2) \times E$. Let $w_1,\,
 w_2$ be coordinates on ${\Bbb C}^2$ and $w_3$ be a coordinate on $E$. Define
 $U=|w_1|^2 + |w_2|^2$ and 

$$ f_a(U)=U
 \sqrt{1+\frac{a^2}{U^2}}+a \ln \frac{U}{\sqrt{U^2+a^2}+a}, \, a>0. \tag 1$$

\noindent Then $\tld{g}_a=\sqrt{-1}\parl\bar{\parl}(f_a(U)+|w_3|^2)$ are Ricci-flat
 K\"{a}hler metrics on $K_{{\Bbb C}P^1} \times E$. Note that we need to extend these
 metrics under blowup. These metrics are asymptotically flat at $\infty$.

Now we construct the approximate Ricci-flat metrics on $M$. We glue $M_0$ 
with 8 copies of $K_{{\Bbb C}P^1} \times E_2$ and  8 copies of
$K_{{\Bbb C}P^1} \times E_3$, by patching the boundaries of some fixed tubular 
neighborhoods of $F$ in $M_0$, the boundaries of some fixed tubular 
neighborhoods of ${\Bbb C}P^1 \times E_2$ in $K_{{\Bbb C}P^1} \times E_2$ and the 
boundaries of some fixed tubular neighborhoods of
${\Bbb C}P^1 \times E_3$ in $K_{{\Bbb C}P^1} \times E_3$. Topologically
the result manifold is $M$. We glue the K\"{a}hler potential of the flat 
orbifold metric on $M_0$ and the 16 copies of the K\"{a}hler potential described 
above to get a function $h_{\vec{ a}}$ (we use different $a$'s for different divisors,
$\vec{a}=(a_1, \cdots, a_{16})$).
Since all metrics are nearly flat near the boundaries as long as $\vec{ a}$ is 
very small, $g_{\vec{a}}=\sqrt{-1} \parl\bar{\parl}h_{\vec{ a}}$ are  K\"{a}hler
metrics on $M$.  They are very close to Ricci flat (as close as we want by choosing
 small $\vec{ a}$). From Yau's existence theorem of Ricci-flat  K\"{a}hler  metrics
 (see [Y]), there is a unique function $u_{\vec{ a}}$ 
on $M$ with $\int_M u_{\vec{ a}}dg_{\vec{ a}} =0$ such that $g^{RF}_{\vec{ a}}$ is
 Ricci flat, where $g^{RF}_{\vec{ a}} =g_{\vec{ a}}+\sqrt{-1} \parl \bar{\parl}
u_{\vec{ a}}$.  

We need the following

\pro{Theorem 2} Let $F$ be the set of singular points in $M_0$. For any relatively 
compact set $W$ in the complement of the proper transformation of $F$ in $M$, 
there exists positive constant $C$ independent of $\vec{a}$ but depending on $W$,
 such that

$$||u_{\vec{ a}}||_{\tld{C}^{4,\alp}(W)} \leq C \cdot |{\vec{ a}}|^2 ,  \tag 2 $$

\noindent where $\tld{C}^{4,\alp}(W)$ is the H\"{o}lder norm with respect to some
fixed coordinate  system on $M$.
\endpro

\proof The approximate Ricci-flat metrics $g_{\vec{ a}}$ on $M$ are similar to the 
approximate Ricci-flat metrics $\om_{a}$ in {\cite L}. The same proof of Theorem 3.4 
in {\cite L} gives the proof of the theorem  here. 
\qed

{\bf 1c}. Special Lagrangian tori on $M$. First 
$E_1 \times E_2 \times E_3$ has special Lagrangian torus fibration with respect to
 holomorphic (3,0)-form $dz_1 \wed dz_2 \wed dz_3$ and the flat metric, 
namely $T_{\alp,\beta,\ga} =T_{\alp} \times T_{\be} \times T_{\ga}$ for any real
 numbers $\alp,\, \be$ and $\ga$, where $T_{\alp} \sst E_1$ is the image of $\alp +
 i{\Bbb R}$ under the projection ${\Bbb C} \rarw E_1$ (Here we need the 
periods $\tau_i$ to be pure imaginary, $i=1,2,3$). For generic values of  $\alp,\, 
\be$ and $\ga$, the image of $T_{\alp,\be,\ga}$ in $M_0$
does not intersect with the singular set $F$. They are embedded
tori in $M_0$. Now we conclude that these tori can be perturbed to special 
 Lagrangian tori in $M$ (embedded).

\pro{Theorem 3} Any special Lagrangian torus $f_0$ in $M_0$ as described above,
 can be perturbed to a  special  Lagrangian torus in $M$.
\endpro

\proof Assume that open set $U$ contains the image of $f_0$ in  $M_0 $,
 the closure $\bar{U}$ is compact. On $\bar{U}$, the metric $g^{RF}_{\vec{ a}}$
 differs from the flat metric on  $E_1\times E_2 \times E_3 $ by an exact form. 
The difference is small on  $\bar{U}$ by Theorem 2. $\Om$ on $\pi^{-1}(\bar{U})$ is 
the same as $\Om_0$ on $\bar{U}$. The image of $f_0$ in $M$ is approximate  special
Lagrangian torus. Now we can apply the proof of Theorem 2.1 in {\cite L} to conclude
that $f_0$ can be perturbed to a  special  Lagrangian torus provided we choose
${\vec{ a}}$ small enough.
\qed

Next we remove the assumption that the periods $\tau_i$ being pure imaginary. 
We can view the threefolds with general $\tau_i$  as deformations of those with pure
imaginary periods. Then by applying Theorem 2.1 i) in {\cite L}, we conclude that the
threefolds have a family of embedded special Lagrangian tori when the real parts of 
$\tau_i$ are sufficiently small.  

We check that the special Lagrangian torus $f: T^3 \raw M$ satisfies 
 $f^*H^2(M)=0$. This is the condition required by mirror symmetry (see
\cite{L}). Note that $h^2(M)=19$ by \cite{B}. The following are a basis 
of $H^2(M)$. There 
are 16 classes which are the Poincare dual of the exceptional divisors
in $M$. Since the image $f(T^3)$ has no 
intersection with the exceptional divisors, the pull-backs of these classes are
zero. Again let
$z_1=x_1+iy_1,\, z_2=x_2+iy_2,\, z_3=x_3+iy_3$ be the complex coordinates 
of $E_1 \times E_2 \times E_3$. Then $dx_1 \wed dy_1, \,dx_2 \wed dy_2, \,dx_3
 \wed dy_3 $ are the classes of degree two invariant under actions 
$\alp$ and $\beta$.
They can be lifted to three classes in $H^2(M)$. Obviously their pull-backs on $T^3$
are zero classes.

\head 2. Special Lagrangian Submanifolds on $K_{{\Bbb C}P^n}$
\endhead

In the following the periods $\tau_1, \, \tau_2 $  and $\tau_3$ are pure imaginary. 
From section 1c we know that one can perturb any special Lagrangian torus $T_{\alp,
\beta,\ga}$ to one in $M$ as long as $\vec{a}$ is small, except $(\alp, \be)=
 (\frac{1}{4},\frac{1}{4}),\, (\frac{1}{4},\frac{3}{4}),\, (\frac{3}{4},\frac{1}{4}),
\,(\frac{3}{4},\frac{3}{4})$ or  $(\alp, \ga)= (0,0),\, (0,\frac{1}{2}),\, 
(\frac{1}{2},0), \,(\frac{1}{2},\frac{1}{2})$ (these are also tori which are preserved 
by either action $\alp$ or $\be$). Note that each such torus intersects 4 fixed elliptic
curves in $F$. For example, $T_{0,\be,0}$ intersects $0 \times E_2 \times 0, \, 0
 \times E_2 \times \frac{\tau_3}{2},\,  \frac{\tau_1}{2} \times E_2 \times 0, \,
\frac{\tau_1}{2}\times E_2 \times \frac{\tau_3}{2}$. Here we try to describe what 
happens to these special Lagrangian tori.

Let $z_1=u_1+ iv_1, \,z_2 =u_2+ iv_2$ be coordinate on ${\Bbb C}^2$.
A obvious family of special Lagrangian submanifolds $L^0_{bc}$ on ${\Bbb C}^2$ are 
given by: $u_1+iu_2 =(b+ic)(v_2+iv_1)$. It is clear that they are invariant under the
${\Bbb Z}_2$ action of ${\Bbb C}^2$. We show 

\pro{Theorem 4} Under the blow up of
 ${\Bbb C}^2/{\Bbb Z}_2$ (the blowup is $K_{{\Bbb C}P^1}$), $L^0_{bc}$ give a family of
special Lagrangian which covers $K_{{\Bbb C}P^1}$. 
\endpro

\proof First we show that there are submanifolds $L_{bc}$ in  $K_{{\Bbb C}P^1}$ 
corresponding to extending $L^0_{bc}$ across the exceptional divisor. Blowup
coordinates are $p=p_1+ ip_2, \, q=q_1+ iq_2$ with  $z_1=p^{1/2},
 \, z_2 =p^{1/2}q$. Combining with the equation of $L^0_{bc}$ we get by eliminating
$u_1, u_2, v_1, v_2$

$$\align
& bq_1^2+bq_2^2 -2bq_1+ (b^2+c^2-1)q_2-b  =0,  \tag 3  \\
& \frac{(bq_1-c)^2-(bq_2-1)^2}{(bq_1-c)^2+(bq_2-1)^2} = \frac{p_1}{\sqrt{p_1^2+p_2^2}}.
 \tag 4
\endalign
$$

\noindent Equation (3) shows that the intersection of special Lagrangian $L_{bc}$ with 
${\Bbb C}P^1$ is a circle. So the topological type of $L_{bc}$ is $S^1 \times
{\Bbb R}$ for all $b$ and $c$. Note that  $L_{bc}$ cover 
${\Bbb C}P^1$ for  $c=0$ and $-1 \leq b \leq 1$. It is not a fibration and there is no
sub-family of $L_{bc}$ to form a fibration of ${\Bbb C}P^1$ with fibre $L_{bc} \cap
{\Bbb C}P^1$.

 Similar to {\bf 1a} the holomorphic (2,0)-form $\Om$ on $K_{{\Bbb C}P^1}$  is 
induced from $dz_1 \wed dz_2$. So Im$\Om|_{L_{bc}}=0$. 

The Ricci-flat metric on $K_{{\Bbb C}P^1}$ is given by

$$\om= \frac{\sqrt{-1}}{U^2 \sqrt{1+U^2}}[(1+U^2)U \parl \bar{\parl}U
-\parl U \wed \bar{\parl}U].$$

\noindent where $U=|z_1|^2+|z_2|^2$. A easy calculation shows that
$\parl \bar{\parl}U|_{L_{bc}}=0$ and $\parl U|_{L_{bc}}$ is real one form. So 
$\om|_{L_{bc}}=0$. $L_{bc}$ are special Lagrangian submanifolds. 
\qed

Note that special Lagrangian $L^{0}_{00} \times T_{\be}$ matches with 
special  Lagrangian 
$T_{0\be 0}$ near each singular locus  $0 \times E_2 \times 0, \, 0
 \times E_2 \times \frac{\tau_3}{2},\,  \frac{\tau_1}{2} \times E_2 \times 0$, and
$\frac{\tau_1}{2}\times E_2 \times \frac{\tau_3}{2}$ in $M_0$. We glue four 
copies of $L_{00} \times T_{\be}$ to $T_{0\be 0}$ and get $\tilde{T}_{0 \be 0}$
in $M$,  it is tempting to
 think that some perturbation of $\tilde{T}_{0 \be 0}$ will be a  special
 Lagrangian torus in $M$. We can not
 prove it because we do not know 
how to estimate the first eigenvalue of Laplace operator acting on $\Om^1(
\tilde{T}_{0 \be 0})$.

Finally we construct special Lagrangian submanifolds in $K_{{\Bbb C}P^{n-1}}$ which is 
isomorphic to the blowup of ${\Bbb C}^{n}/{\Bbb Z}_{n}$. Let $z_1=x_1+iy_1, 
\cdots, z_{n}=x_{n}+iy_{n}$ be complex coordinates of  ${\Bbb C}^{n}$. The 
holomorphic $(n,0)$-form $\Om$ on $K_{{\Bbb C}P^{n-1}}$ is the pull-back of $dz_1
\wed \cdots \wed dz_n$. Let $U=|z_1|^2 + \cdots + |z_n|^2$. The Ricci-flat
 K\"{a}ler form on 
$K_{{\Bbb C}P^{n-1}}$ is

$$\align
&\om=\sqrt{-1}(f'(U) \parl \bar{\parl} U + f''(U) \parl U
\wed  \bar{\parl} U), \\
&f'(U)= (1+\frac{1}{U^n})^{1/n}.
\endalign
$$

Let $\vec{x}=(x_1,\cdots,x_n),\, \vec{y}=(y_1,\cdots,y_n)$ and $A=(a_{ij})_{n \times n}$
a real matrix. Consider the real $n$-dimensional plane $L^0_{A}$ in ${\Bbb C}^{n}$ 
defined by $\vec{x}=A \vec{y}$. It is easy to check that $ \parl \bar{\parl} 
U|_{L^0_{A}}=0$ and $\parl U|_{L^0_A}$ is real one form if and only if $A$ 
equals its transpose $A^t$.
So $A=A^t$ implies that $\om|_{L^0_A}=0$. Note that the only
singular point on the image of 
$L_A^0$ in ${\Bbb C}^{n}/{\Bbb Z}_{n}$ is origin.

Let $z_1=w_1^{1/n}, z_2=w_1^{1/n}w_2, \cdots, z_n=w_1^{1/n}w_n$ be the blowup 
coordinates. Let $L_A$ be the image of $L^0_A$ under blowup. We show that 
$L_A$ is smooth.
Rewrite $\vec{x}= A \vec{y}$ in terms of $z_1=x_1+iy_1, w_2=u_2+iv_2, \cdots, w_n=
u_n +iv_n$, we have

$$\align
&(1-\sum_{j=2}^n a_{1j}v_j)x_1 = (a_{11}+\sum_{j=2}^n a_{1j}u_j)y_1, \tag 5 \\
&(u_i-\sum_{j=2}^n a_{ij}v_j)x_1=(a_{i1}+v_i+\sum_{j=2}^n a_{ij}u_j)y_1, \,\,\, 2 \leq
i \leq n. \tag 6i
\endalign
$$  

\noindent Dividing (6i) by (5) we get $n-1$ equations defining the intersection of $L_A$
with ${\Bbb C}P^{n-1}$ which is smooth. In particular when $n=3$, $a_{11}=a_{22}=1$ 
and all other $a_{ij}=0$, we get the equations $u_2+v_2=0, \, u_3+v_3=0$ which defines 
a $S^1 \times S^1$ in ${\Bbb C}P^{2}$. So the topological type of $L_A$ is $S^1 \times
S^1 \times {\Bbb R}$ for the $A$.

Let $A_{i_1 \cdots i_k}$ be a $k \times k$ matrix formed from elements in $A$ with both
row and column numbers in set $\{i_1, \cdots,i_k \}$. Then 

$$
dz_1\wed \cdots \wed dz_n |_{L^0_A} = (\sum_{k=\text{even}} 
\sum_{i_1 < \cdots <i_k}     
i^k \det (A_{i_1 \cdots i_k}) + \sum_{k=\text{odd}} \sum_{i_1 < \cdots <i_k}     
i^k \det (A_{i_1 \cdots i_k}))dy_1 \wed \cdots \wed dy_n.
$$

Taking the phase factor into consideration we conclude that $L_A$ is special 
Lagrangian submanifolds in $K_{{\Bbb C}P^{n-1}}$ when $A= A^t$ and

$$\sin \hta \cdot  \sum_{k=\text{even}} 
\sum_{i_1 < \cdots <i_k} \det (A_{i_1 \cdots i_k}) + \cos \hta
 \sum_{k=\text{odd}} \sum_{i_1 < \cdots <i_k} \det (A_{i_1 \cdots i_k}), $$ 

\noindent where $0 \leq \hta < 2 \pi$. Furthermore any real $n$-dimensional
hyperplane in ${\Bbb C}^n$, which is 
the limit of special Lagrangian $L^0_{A}$, gives rise to a special Lagrangian 
submanifold in $K_{{\Bbb C}P^{n-1}}$.


\Refs

\widestnumber\key{FMS$^2$}

\ref\key B \by C. Borcea \paper
  K3 surfaces with involution and mirror pairs of Calabi-Yau manifolds
\jour in Mirror Symmetry II, ed. by B. Greene and S.T. Yau 
\publ International Press
\pages 717-743   \yr 1997 
\endref

\ref\key Br \by R. Bryant \paper
Some Examples of Special Lagrangian Tori
\jour preprint   \yr 1999 
\endref

\ref\key GW \by M. Gross and P.M.H. Wilson \paper
Mirror symmetry via 3-tori for a class of Calabi-Yau
threefolds 
\jour Math. Ann. \vol 309 \pages 505-531 \yr 1997
\endref


\ref\key L \by  P. Lu \paper
K\"{a}hler-Einstein metrics on Kummer threefold and 
special Lagrangian tori
\jour preprint, to appear on Comm. Anal. Geom. \yr 1997
\endref


\ref\key SYZ \by A. Strominger, S.T. Yau and E. Zaslow \paper
Mirror Symmetry is T-Duality 
\jour Nucl. Phys. \vol B479 \pages 243-259  \yr 1996
\endref



\ref\key V \by C. Voisin \paper
Miroirs et involutions sur les surfaces K3   
\jour Ast\'{e}risque
\vol 218 \pages 273-323 \yr 1993
\endref


\ref\key Y \by S.T. Yau \paper
On the Ricci Curvature of a compact K\"{a}hler manifold and the complex 
Monge-Amp\'{e}re equations. I 
\jour Comm. Pure Appl. Math. 
\vol 31 \pages 339-411 \yr 1978  
\endref
\end{document}